\documentclass[12TP,draft]{article}

\begin{document}

\begin{center}
\LARGE\noindent\textbf{A sufficient condition for the hamiltonian property of digraphs with large semi-degrees}\\

\end{center}
\begin{center}
\noindent\textbf{S. Kh.                                                                                                                                                                                                                                                                                                                                                                                                     Darbinyan }\\

Institute for Informatics and Automation Problems\\ Armenian National Academy of Sciences, P. Sevak 1, Yerevan 0014, Armenia

e-mail: samdarbin@ipia.sci.am\\
\end{center}

\textbf{Abstract}\\

 Let $D$ be a digraph on $p\geq 5$ vertices with  minimum degree at least $p-1$ and with minimum semi-degree at least $p/2-1$. For $D$ (unless  some extremal cases) we present a detailed proof of the following results [12]: (i) $D$ contains cycles of length 3, 4 and $p-1$; (ii) if $p=2n$, then $D$ is hamiltonian. \\


\noindent\textbf{1. Introduction and Terminology}\\

Ghouila-Houri [18] proved that every strong digraph on $p$ vertices and with minimum degree at least $p$ is hamiltonian. There are many extentions of this theorem for digraphs and orgraphs. In particular, in many papers, various degree conditions have been obtained for digraphs (orgraphs) to be hamiltonian or pancyclic or vertex pancyclic (see e.g. [2]-[33]). C. Thomassen [31] proved that  any digraph on $p=2n+1$ vertices with minimum semi-degree at least $n$ is hamiltonian unless  some extremal cases, which are characterized. In [9], we proved that if a digraph $D$ on $2n+1$ vertices satisfies the conditions of this Tomassen's theorem, then $D$ also is pancyclic (the extremal cases are characterized). For additional information on hamiltonian and pancyclic digraphs, see the book [1] by J. Bang-Jenssen and G. Gutin.

 In this paper we present a detailed proof of the following results. 

 Every digraph $D$ (unless some extremal cases) on  $p\geq 5$ vertices with minimum degree at least $p-1$ and with minimum semi-degree at least $p/2-1$: (i) $D$ has  cycles of length 3, 4 and $p-1$; (ii) if $p=2n$, then $D$ is hamiltonian (in [12], we gave only a short outline of the proofs of this results). 
 
In this paper we will consider finite digraphs without loops and multiple arcs. We denote the vertex set of digraph $D$ by $V(D)$ and its arc set by  $A(D)$. We will often use $D$ instead of $A(D)$ and $V(D)$. The arc from a vertex $x$ to a vertex $y$ will be denoted by $xy$. If $xy$ is an arc, then we say that $x$ dominates $y$ (or $y$ is dominated by $x$). For $A$,  $B\subset V(D)$, we define $A(A\rightarrow B)$  as the set $\{xy\in A(D) / x\in A, y\in B\}$ and $A(A,B)=A(A\rightarrow B)\cup A(B\rightarrow A)$. If  $x\in V(D)$ and $A=\{x\}$, we often write $x$ instead of $\{x\}$. For disjoint subsets  $A$ and $B$  of $V(D)$,  $A\rightarrow B$ means that every vertex of $A$ dominates every vertex of $B$. If $C\subset V(D)$, $A\rightarrow B$ and $B\rightarrow C$, then we write $A\rightarrow B\rightarrow C$. The outset of vertex $x$ is the set $O(x)=\{y\in V(D) / xy\in A(D)\}$ and $I(x)=\{y\in V(D) / yx\in A(D)\}$ is the inset of $x$. Similarly, if $A\subseteq V(D)$ then $O(x,A)=\{y\in A / xy\in A(D)\}$ and $I(x,A)=\{y\in A / yx\in A(D)\}$. The out-degree of $x$ is $od(x)=|O(x)|$ and $id(x)=|I(x)|$ is the in-degree of $x$. Similarly, $od(x,A)=|O(x,A)|$ and $id(x,A)=|I(x,A)|$. The degree of the vertex $x$ in $D$ is defined as $d(x)=id(x)+od(x)$. The subdigraph of $D$ induced by a subset $A$ of $V(D)$ is denoted by $\langle A\rangle$. All paths and cycles we consider in this paper are directed and simple. The path ( respectively, the cycle ) consisting of distinct vertices $x_1,x_2,\ldots ,x_n$ ( $n\geq 2 $) and arcs $x_ix_{i+1}$, $i\in [1,n-1]$  ( respectively, $x_ix_{i+1}$, $i\in [1,n-1]$, and $x_nx_1$ ), is denoted by $x_1x_2\ldots x_n$ (respectively, $x_1x_2\ldots x_nx_1$ ). The cycle on $k$ vertices is denoted $C_k$. For a cycle  $C_k=x_1x_2\ldots x_kx_1$, we take the indices modulo $k$, i.e., $x_s=x_i$ for every $s$ and $i$ such that $i\equiv s \,\hbox {mod}\, k$.

  Two distinct vertices $x$ and $y$ are adjacent if $xy\in A(D)$ or $yx\in A(D) $ (or both) (i.e, $x$ is adjacent with $y$ and $y$ is adjacent with $x$). The notation $A(x,y)\not= \emptyset$ (respectively, $A(x,y)=\emptyset $) means that the vertices $x$ and $y$ are adjacent (respectively, are not adjacent).

 The converse digraph $\overleftarrow {D}$ of a digraph $D$ is the digraph obtained from $D$ by reversing all arcs of $D$.

 For an undirected graph $G$, we denote by $G^*$ symmetric digraph obtained from $G$ by replacing every edge $xy$ with the pair $xy$, $yx$ of arcs. Further, $C^*(5)$ is a symmetric digraph  obtained from undirected cycle of length 5. $K_n$ (respectively, $K_{n,m}$)  denotes the complete undirected graph on $n$ vertices (respectively, undirected complete bipartite graph, with partite sets of cardinalities $n$ and $m$), and  $\overline K_n$ denotes the  complement of $K_n$.   If $G_1$ and $G_2$ are undirected graphs, then $G_1\cup G_2$ is the disjoint union of $G_1$ and $G_2$. The join of $G_1$ and $G_2$, denoted by  $G_1 + G_2$, is the  union of $G_1\cup G_2$ and of all the edges between $G_1$ and $G_2$.

For integers $a$ and $b$, let $[a,b]$  denote the set of all integers which are not less than $a$ and are not greater than $b$.
 We refer the reader to J.Bang-Jensens and G.Gutin's book [1] for notations and terminology not defined here. \\

\noindent\textbf{2. Preliminaries and Additional notations}\\

Let us recall some well-known lemmas used in this paper.\\

\noindent\textbf{Lemma 1} ([21]). Let $D$ be a digraph on $p\geq 3$ vertices containing a cycle $C_m$, $m\in [2,p-1] $. Let $x$ be a vertex not contained in this cycle. If $d(x,C_m)\geq m+1$, then  $D$ contains a cycle $C_k$ for all  $k\in [2,m+1]$. \framebox  \\\\

The following Lemma will be used extensively in the proofs our results.\\

\noindent\textbf{Lemma 2} ([6]). Let $D$ be a digraph on $p\geq 3$ vertices containing a path $P:=x_1x_2\ldots x_m$, $m\in [2,p-1]$ and let $x$ be a vertex not contained in this path. If one of the following holds:

 (i) $d(x,P)\geq m+2$; 

 (ii) $d(x,P)\geq m+1$ and $xx_1\notin D$ or $x_mx_1\notin D$; 

 (iii) $d(x,P)\geq m$, $xx_1\notin D$ and $x_mx\notin D$;

\noindent\textbf{}then there is an  $i\in [1,m-1]$ such that $x_ix,xx_{i+1}\in D$, i.e., $D$ contains a path $x_1x_2\ldots x_ixx_{i+1}\ldots x_m$ of length $m$  (we say that  $x$ can be inserted into $P$ or the path  $x_1x_2\ldots x_ixx_{i+1}\ldots x_m$ is extended from $P$ with  $x$). \fbox \\\\

As an immediate consequence of Lemma 2, we get the following:\\

\noindent\textbf{Lemma 3}. Let $D$ be a digraph on $p\geq 4$ vertices and let  $P:=x_1x_2\ldots x_m$, $m\in [2,p-2]$, be a path of  maximal length from $x_1$ to $x_m$ in $D$. If the induced subdigraph $\langle V(D)\setminus V(P)\rangle$ is strong and  $d(x,V(P))=m+1$ for every vertex $x\in V(D)\setminus V(P)$, then there is an integer $l\in [1, m]$ such that $O(x,V(P))=\{x_1,x_2,\ldots ,x_l\}$ and $I(x,V(P))=\{x_l,x_{l+1},\ldots , x_m\}$.  \framebox \\\\

Now we introduce the following notations.\\

 \noindent\textbf {Notation}. For any positive integer $n$, let $H(n,n)$ denote the set of digraphs $D$ on  $2n$ vertices  such that 
$V(D)=A\cup B$,\, $\langle A \rangle \equiv \langle B \rangle \equiv K_{n}^*$, \, $A(B\rightarrow A)=\emptyset$\, and for every vertex $x\in A$ (respectively,\, $y\in B$) \,$A(x\rightarrow B)\not=\emptyset$ (respectively,\, $A(A\rightarrow y)\not=\emptyset$).\\

\noindent\textbf {Notation}. For any integer $n\geq 2$, let $H(n,n-1,1)$ denote the set of digraphs $D$ on $2n$ vertices such that \, $V(D)= A\cup B\cup \{a\}$\,, $|A|=|B|+1=n$,\, $A(\langle A\rangle )=\emptyset$, \,$\langle B\cup \{a\}\rangle \subseteq K_n^* $, $yz,\,zy\in D$  for each pair of vertices $y\in A$,  $z\in B$  and either $I(a)=B$ and $a\rightarrow A$ or $O(a)=B$ and $A\rightarrow a$.\\

\noindent\textbf {Notation}. For any integer $n\geq 2$ define the digraph $H(2n)$ as follows: \,$V(H(2n))=A\cup B\cup \{x,y\}$, \, $\langle A\rangle \equiv \langle B\rangle \equiv K_{n-1}^*$,\, $A(A,B)=\emptyset$, \, $O(x)=\{y\}\cup A$,\, $I(x)=O(y)=A\cup B$\, and \, $I(y)=\{x\}\cup B$. 

$H'(2n)$ is a digraph obtained from $H(2n)$ by adding the arc $yx$.\\

\noindent\textbf{Notation.} Let $D_6$ be a digraph with vertex set $\{x_1,x_2,\ldots ,x_5,x\}$  and arc set 
$$
\{x_ix_{i+1} \,/ 1\leq i\leq 4 \}\cup \{xx_i/\,1\leq i\leq 3\} \cup \{x_1x_5,x_2x_5,x_5x_1,x_5x_4,x_3x_2,x_3x,x_4x_1,x_4x \}.$$ By $D'_6$ we denote a digraph obtained from $D_6$ by adding the arc  $x_2x_4$.\\

 Note that the digraphs  $D_6$ and $D'_6$ both are not hamiltonian and each  of $D_6$ and $D'_6$ contains a cycle of length 5. \\ 

\noindent\textbf{Lemma 4}. Let $D$ be a digraph on $p\geq 3$ vertices   with  minimum degree  at least $p-1$ and with minimum semi-degree  at least $p/2-1$. Then

(i) either $D$ is strong or $D\in H(n,n)$;

(ii) if $B\subset V(D)$,  $|B|\geq (p+1)/2$ and  $x\in V(D)\setminus B$, then  $A(x\rightarrow B)\not=\emptyset$ and $A(B\rightarrow x)\not=\emptyset$. \framebox \\\\

 \noindent\textbf {3. A sufficient condition for the existence of cycles of length $|V(D)|-1$ in digraph $D$}\\
   
\noindent\textbf{Theorem 1}. Let $D$ be a digraph on $ p\ge5$ vertices  with minimum degree  at least $p-1$ and with minimum semi-degree  at least $p/2-1$. Then $D$ has a cycle  of length $p-1$ unless 
$$
D\in H(n,n)\cup \{[(K_n\cup K_n)+K_1]^*, H(2n), H^\prime (2n), C^*(5)\} \quad  \hbox {or else} \quad p=2n \quad \hbox {and} \quad D\subseteq K^*_{n,n}.
$$
    
\noindent\textbf {Proof}. By Lemma 4(i), the result is easily verified if $D$ is not strong. Assume that $D$ is strong. Suppose, on the contrary, that the theorem is not true. In particular, $D$  contains no cycle of length $p-1$. Let $ C:= C_{m}:= x_1x_2\ldots x_mx_1$ be an arbitrary non-hamiltonian cycle of maximum length in $D$. It is easy to see that $m\in [3,p-2]$.

  From Lemma 1 and the maximality of $m$ it follows that for each vertex $y\in B:=V(D)\setminus V(C)$ and for each  $i\in [1,m]$,
$$
 d(y,C)\leq m,  \quad  d(y,B)\geq p-m-1  \quad \hbox {and if} \quad x_iy\in D, \quad \hbox {then } \quad yx_{i+1}\notin D.
 \eqno{(1)}
$$
 Using $d(y,B)\geq p-m-1$ it is not difficult to show the following claim:\\

\noindent\textbf {Claim 1}. For any two distinct vertices $x, y\in B$ if in subdigraph $\langle B \rangle $ there is no path from  $x$ to $y$, then in $\langle B \rangle $  there is a path from $y$ to $x$ of length at most 2. \fbox \\\\

We first prove the following Claims 2 and 3.\\ 

\noindent\textbf {Claim 2}. The induced subdigraph $\langle B \rangle $ is strongly connected.

\noindent\textbf {Proof}. Suppose, on the contrary, that $\langle B \rangle $ is not strong. Let $D_1$, $D_2$, \ldots  ,     $D_s$ $(s\geq 2)$ be the strong components of $\langle B \rangle $  labeled in such  a way that no vertex of $D_i$ dominates a vertex of $D_j$ whenever $i>j$. From Claim 1 it follows  that for each pair of vertices $ y\in V(D_1)$ and $z\in V(D_s)$ in  $\langle B \rangle $ there is a path  from  $y$ to  $z$ of length 1 or 2. We choose the vertices $ y\in V(D_1)$ and $z\in V(D_s)$ such that the path $y_1y_2\ldots y_k$, where $y_1:=y$ and $y_k:=z$, will have minimum length  among all paths in $\langle B \rangle $  with origin vertex in $D_1$ and terminus vertex in $D_s$. By Claim 1, $k=2$ or $k=3$. We consider the following tree cases.\\

\noindent\textbf {Case 1}. $k<|B|=p-m.$

It follows from the maximality of  $C$ that if $x_iy_1\in D$, where  $i\in [1,m]$, then $ A(y_k\rightarrow \{ x_{i+1},x_{i+2},$  $\ldots ,x_{i+k}\})=\emptyset$. Since $D$ is strong, we see that $C\not\subseteq I(y_1)$. Therefore the vertex $y_k$ dose not dominate at least $ id(y_1,C)+1$ vertices of $C$. On the other hand, we have $A(y_k\rightarrow V(D_1))=\emptyset$ and $ I(y_1)\subset C\cup V(D_1)$. Hence the vertex $y_k$ dose not dominate at least $id(y_1)+3$ vertices. From this we obtain $od(y_k)\leq p-id(y_1)-3\leq p/2-2$, which is a contradiction.\\ 

\noindent\textbf {Case 2}. $k=|B|=2$.

It is easy to see that $s=2$, $m=p-2$, $V(D_1)=\{ y_1 \}$, $V(D_2)=\{ y_2 \}$, $ I(y_1)\subset C$ and 
$$
| A(x_i\rightarrow y_1)| +  |A(y_2\rightarrow x_{i+2})| \leq 1      
$$
for all $i\in [1,m]$. Hence the vertex $y_2$ dose not dominate at least $id(y_1)+2$ vertices. Therefore $od(y_2)\leq p-id(y_1)-2\leq p/2-1$. It follows that $p=2n$, $id(y_1)=od(y_2)=n-1$  and
$$
y_2x_i\in D \quad \hbox {if and only if} \quad  x_{i-2}y_1 \notin D.      \eqno {(2)} 
$$
By Lemma 1, it is easy to see that $d(y_1)=d(y_2)=2n-1$ and $od(y_1)=id(y_2)=n$, $m\geq 4$. Now we divide this case into two subcases. \\

\noindent\textbf {Subcase 2.1}.  $y_1\rightarrow \{ x_i,x_{i+1} \}$ for some $i\in [1,m]$.

Note that, by Lemma 2, without loss of generality,  we may assume  that $x_my_1 \in D, $ $y_1\rightarrow \{x_2,x_3\}$  and  $A(x_1,y_1)=\emptyset .$ From this, (1) and (2) it follows that $x_2y_1 \notin D$, $y_2x_3$, $y_2x_4\in D$ and $A(x_2,y_2)=\emptyset$. Therefore, by Lemma 2 we have $x_1y_2\in D$ since $d(y_2,C)=2n-2$ and the vertex $y_2$ cannot be inserted into the path $x_3x_4 \ldots x_mx_1$. If $x_2x_1 \in D$, then $ C_{2n-1}=x_my_1x_2x_1y_2x_4\ldots x_m$. This contradicts our supposition that $D$  contains no cycle of length $p-1$. Hence, $x_2x_1\notin D$. From this and $A(x_2,y_2)=A(x_2\rightarrow y_1)=\emptyset$ it follows that $d(x_2,\{x_3,x_4,\ldots ,x_m \})\geq 2n-3$. Therefore by Lemma 2, $x_mx_2 \in D$ since the vertex $x_2$ cannot be inserted into the path $x_3x_4\ldots x_m$. Now it is easy to see that $| A(x_i\rightarrow y_1)| + | A(y_2\rightarrow x_{i+1})| \leq 1$ for all $i\in [2,m-1]$. 
Therefore $x_3y_1\notin D$ since $y_2x_4\in D$. From this and (2) it follows that $y_2x_5\in D$ and $x_4y_1\notin D$. Continuing in this manner, we obtain that $A(\{x_5,x_6,\ldots, x_{m-1}\}\rightarrow y_1)=\emptyset $. Therefore  $A(\{x_1,x_2,\ldots, x_{m-1}\}\rightarrow y_1)=\emptyset $, which is a contradiction.\\

\noindent\textbf {Subcase 2.2}. $|A(y_1\rightarrow \{ x_i,x_{i+1})| \leq 1$ for all $i\in [1,m]$.

Since $od(y_1)=n$,  we can assume that  $O(y_1)=\{x_1,x_3,\ldots ,x_{2n-3},y_2 \}$. Using this and  $od(y_2)=id(y_1)=n-1$, we obtain $I(y_1)=\{x_1,x_3 ,\ldots , x_{2n-3}\}$. Therefore by (2),
 $$
O(y_2)=\{x_2,x_4,\ldots ,x_{2n-2}\} \quad \hbox {and} \quad  I(y_2)=\{y_1,x_2,x_4,\ldots ,x_{2n-2}\}.
$$
If $x_ix_j\in D$ for distinct vertices   $x_i,x_j \in \{x_1,x_3,\ldots ,x_{2n-3}\}$, then  $C_{2n-1}=y_1x_ix_jx_{j+1}\ldots x_{i-1}y_2x_{i+1}\ldots $ $x_{j-2}y_1$, when  $| \{x_{i+1},x_{i+2}, \ldots ,x_{j-1}\}|  \geq 2$  and $C_{2n-1}=x_ix_jy_1y_2x_{j+1} x_{j+2}\ldots x_{i-1}x_i$,  when $|\{x_{i+1},x_{i+2},$ $\ldots,$ $x_{j-1}\}|  =1$. This contradicts that  $C_{p-1}\not\subset D$. Thus we have
 $$A(\langle \{x_1,x_3,\ldots ,x_{2n-3},y_2 \} \rangle)=\emptyset.$$

Considering the digraph $\overleftarrow {D}$, by the same arguments we obtain  
$$
A(\langle \{x_2,x_4,\ldots ,x_{2n-2},y_1 \} \rangle)=\emptyset.
$$ 
Therefore $D\subseteq K_{n,n}^*$, which contradicts our supposition that the theorem is not true.\\

\noindent\textbf {Case 3}. $k=\mid B\mid =3$.

From the minimality of  $k$ it follows that $y_1y_3\notin D$, $s=3,$  $A(\{y_2,y_3\} \rightarrow  y_1)=\emptyset $ and $V(D_1)=\{ y_1\}$. Hence $I(y_1)\subset C$. On the other hand, from the maximality of the cycle $C$ it follows that for each $i\in [1,m]$  
$$
\hbox {if} \quad x_iy_1\in  D, \quad \hbox{then} \quad A(y_2\rightarrow \{x_{i+1},x_{i+2}\})=\emptyset .
$$
Therefore $y_2$ dose not dominate at least $id(y_1)+3$ vertices, a contradiction.  Claim 2 is proved. \fbox \\\\

\noindent\textbf {Claim 3}. At least two distinct vertices of $C$ are adjacent with some vertices of $B$.

\noindent\textbf {Proof}. Assume that Claim 3 is not true. Then exactly one vertex, say $x$, of  $C$ is adjacent with some vertices of $B$. Hence for each vertex $x_i \in C\setminus \{x\}$ and for each vertex $y\in B$ we have 
 $$
d(x_i)=d(x_i,C)\leq 2m-2 \quad \hbox {and} \quad d(y)=d(y,B)+d(y,x)\leq 2p-2m.
$$
 Since $d(x_i)+d(y)\geq 2p-2$, we conclude that the inequalities above are equalities. This implies that the subdigraphs $\langle C \rangle$ and $\langle B\cup \{x\}\rangle$ are complete. From $d(x_i)=2m-2\geq p-1$ and $d(y)=2p-2m\geq p-1$, we obtain that $p=2m-1$. Therefore $G\equiv [(K_{m-1}\cup K_{m-1})+K_1]^*$, which contradicts our supposition. This  proves Claim 3. \fbox \\\\

Since $D$ is strong, then $A(C\rightarrow B)\not=\emptyset$ and $A(B\rightarrow C)\not=\emptyset$. Together with Claim 3 this  implies  that there are vertices $x_a\not= x_b$, $x_a,x_b\in C$  and  $x,y\in B$ such that $x_ax$, $yx_b\in D$ and
$$
A(\{x_{a+1},x_{a+2},\ldots ,x_{b-1}\},B)=\emptyset, \quad \hbox {if} \quad x_b\not=x_{a+1}.     \eqno {(3)} 
$$
To be definite, assume that $x_b:=x_1$ and $x_a:=x_{m-h}$ ($0\leq h\leq m-2$). We  consider the following two cases.\\

\noindent\textbf {Case 1}. $x_{m-h+1}\not=x_1$ (i.e., $h\geq 1$).

 Consider the paths $P_0$, $P_1$, \ldots , $P_k$ ($0\leq k\leq h$ and $k$ is as maximum as possible), where $P:=P_0:= x_1x_2\ldots x_{m-h}$ and the path $P_i$, $i\in [1,k]$, is extended from the path $P_{i-1}$ with a vertex $z_i\in \{x_{m-h+1}, x_{m-h+2},\ldots , x_m \}\setminus \{z_1,z_2, \ldots, z_{i-1} \}$. Note that the path $P_i$, $i\in [0,k]$, contains $m-h+i$ vertices. It follows that some vertices $y_1,y_2,\ldots ,y_d \in \{x_{m-h+1},x_{m-h+2},\ldots ,x_{m}\}$, where $1\leq d\leq h $, dose not containing the extended path $P_k$. Therefore, using (3) and Lemma 2, for each $z\in B$ and for each $y_i$ we obtain
$$
d(z)=d(z,B)+d(z,C)\leq 2p-2m-2+m-h+1=2p-m-h-1
$$
and
$$
d(y_i)=d(y_i,C)\leq m+d-1.
$$
Hence it is clear that $$2p-2\leq d(z)+d(y_i)\leq 2p+d-h-2.$$

  It is not difficult to see that $h=d$, $d(z,C)=m-h+1$, $d(y_i,C)=m+h-1 $ and the subdigraphs  $ \langle B \rangle $ and $ \langle \{x_{m-h+1},x_{m-h+2},\ldots ,x_{m}\}\rangle $  are complete. By Lemma 2(ii), we also have $ x_{m-h} \rightarrow B \cup \{x_{m-h+1},x_{m-h+2},\ldots ,x_{m}\} \rightarrow x_1 $. It is easy to see that $h=|B|=p-m \geq 2$ and the path $P=x_1x_2\ldots x_{m-h}$ has maximum length among all paths from $x_1$ to $x_{m-h}$ in subdigraph $\langle C \rangle$ and in subdigraph $\langle B \cup \{x_1 ,x_2  ,\ldots ,x_{m-h}\}\rangle $. Therefore by Lemma 3, there are integers $l\in [1,m-h]$ and $ r\in [1,m-h] $ such that 
$$
O(u,P)=\{x_1,x_2,\ldots ,x_l\}, \quad I(u,P)=\{x_l,x_{l+1},\ldots ,x_{m-h}\},
$$
$$
O(z,P)=\{x_1,x_2,\ldots ,x_r\}, \quad I(z,P)=\{x_r,x_{r+1},\ldots ,x_{m-h}\}.\eqno {(4)}
$$
for all $u\in B$ and for all $z\in \{x_{m-h+1},x_{m-h+2},\ldots ,x_{m}\}$.

Without loss of generality, we may assume  that $l\leq r$ (otherwise we will consider the digraph $\overleftarrow {D}$).

 Let $l=1$. Then from $od(u)\geq p/2-1$ and (4) it follows that $h\geq p/2-1$ and $ p\geq  2(p/2-1)+m-h=p-2+m-h$. Since $m-h\geq  2$, we see that $p=2n$, $m-h=2$,\, $h=n-1$ and $r=2$. Therefore $G\in \{H(2n),H '(2n)\}$, which contradicts the our  supposition. 

 Let now $l\geq 2$. We can assume that $r\leq  m-h-1$ (otherwise in digraph $\overleftarrow {D}$ we will have the considered case $l=1$). Since $ \langle \{x_{m-h+1},x_{m-h+2},\ldots ,x_{m}\}\rangle $  are complete and (4), for each vertex  $ z\in \{x_{m-h+1},x_{m-h+2},\ldots ,x_m \}$  we have $I(z)=\{x_r,x_{r+1},\ldots ,x_m \}\setminus \{ z\}$. This implies that $m-r\geq p/2-1$. If $i\in [r+1,m-h]$  and $x_1x_i\in D$ then by (4) and $2\leq l\leq r\leq m-h-1$ we have $ C_{m+1}=x_1x_ix_{i+1}\ldots x_mx_2\ldots x_{i-1}xx_1$, where $x\in B$, a contradiction. Because of this and $2\leq l\leq r$, we may assume that
$$
A(x_1\rightarrow B\cup \{ x_{r+1},x_{r+2},\ldots , x_m \})=\emptyset.
$$
Therefore, since $m-r\geq p/2-1$ and $|B|=h\geq 2$, we obtain $od(x_1)\leq  p-1-h-(m-r)\leq p/2-h$, which contradicts the condition that $od(x_1)\geq p/2-1$.\\

\noindent\textbf {Case 2}. $x_{m-h+1}=x_1$ (i.e., $h=0$).

Then any path from  $x$ to  $y$ in   $\langle B \rangle$ is a hamiltonian path. Let $ u_1u_2\ldots u_{p-m}$ be a hamiltonian path in $\langle B \rangle$, where $u_1:=x$, $u_{p-m}:=y$. From this, if $1\leq i<j\leq p-m$, then  $u_iu_j\in D$ if and only if $j=i+1$.

For this case ($h=0$) we first prove Claims 4-9.\\

\noindent\textbf {Claim 4}. $p-m=2$ (i.e., $m=p-2$).

\noindent\textbf {Proof}. Suppose, to the contrary, that $p-m\geq 3$. It follows from observations above that $u_1u_{p-m}\notin D$ and $od(u_1,B)=id(u_{p-m},B)=1$. From this and (1), we obtain
$$
p-1\leq d(u_1)\leq m+1+id(u_1,B) \quad  \hbox {and} \quad p-1\leq d(u_{p-m})\leq m+1+od(u_{p-m},B).
$$
This implies that $id(u_1,B)$  and  $od(u_{p-m},B)\geq p-m-2$. Therefore in  $ \langle B \rangle$ there is a path from  $u_{p-m}$ to  $u_1$ of length $k=1$ or $k=2$ since $p-m\geq 3$. For any integer $l\geq 1$ , put
$$
I^+_l:=\{x_j/\,x_{j-l}u_{p-m}\in D \}.
$$
Since  $id(u_{p-m},C)=id(u_{p-m})-1$ and \, $C\not\subseteq I(u_{p-m})$, we see that for each $l\in [1,2]$, 
$$
| I^+_l\cup I^+_{l+1} | \geq id(u_ {p-m}).
$$
From the maximality of the cycle $C$ it follows that
$
A(u_1\rightarrow I^+_k\cup I^+_{k+1})=\emptyset.
$
Together with  $A(u_1\rightarrow \{u_3,u_4,\ldots ,$ $u_{p-m}\})=\emptyset$ this implies that 
$$
p/2-1\leq od(u_1)\leq p-1-| I^+_k \cup I^+_{k+1}| -(p-m-2)\leq m+1-id(u_{p-m})\leq m+1-p/2+1.
$$
Therefore, since $m\leq p-3$, we obtain  that $p-m=3$, $p=2n$ and $od(u_1)=id(u_3)=n-1$. Hence,  $id(u_1)$  and $od(u_3)\geq n$. 
We now claim that $u_3u_1$ and $u_2u_1\in D$. Indeed, otherwise $id(u_1,C)\geq n-1$ and if $x_iu_1\in D$, then $A(u_3\rightarrow \{x_{i+2},x_{i+3}\})=\emptyset$. From this it is not difficult to see that  $od(u_3)\leq n-1$, which contradicts the fact that   $od(u_3)\geq n$.

Similarly, we can see that $u_3u_2\in D$. So we have $u_3u_1$, $u_2u_1$, $u_3u_2 \in D$. Then, since $id(u_3,C)=n-2$, $m\geq n$,  $m\geq id(u_3,C)+2$ and $C$ is a  non-hamiltonian cycle of maximal length, it follows that $|\cup _{i=1}^3 I^+_i | \geq n $ and $A(u_1\rightarrow \cup _{i=1}^3 I^+_i)=\emptyset$. Together with $u_1u_3\notin D$ this implies that  $od(u_1)\leq n-2$, a contradiction. This completes the proof of Claim 4. \fbox \\\\

 Note that, by Claims 4 and 2 we have $m=p-2$, $B:=\{u,v\}$ and $uv$, $vu \in D$.\\

\noindent\textbf {Remark}. By symmetry of the vertices $u$ and $v$, Claims 5-9 are also true for the vertex $v$.\\

\noindent\textbf {Claim 5}. If $x_iu$, $ux_{i+2} \in D$, $i\in [1,m]$, then $ |A(x_{i+1},v)| =2$ (i.e., $x_{i+1}v$ and $vx_{i+1}\in D$).

\noindent\textbf {Proof}. Since the cycle $x_iux_{i+2}x_{i+3}\ldots x_i$ has length $m$ and the vertices $v$ and $x_{i+1}$ are not on this cycle, the subdigraph $\langle \{v,x_{i+1} \} \rangle$ is strong  by Claim 2. Therefore $vx_{i+1}$ and $x_{i+1}v \in D$. \fbox \\\\

From Claim 5, $uv$, $vu\in D$ and the maximality of the cycle $C$ we have the following:\\

\noindent\textbf {Claim 6}. If $i\in [1,m]$, then
$$
| A(\{x_i,x_{i+1}\}\rightarrow u)| +| A(u\rightarrow x_{i+3})| \leq 2 \quad \hbox {and}\quad  | A(x_{i-2}\rightarrow u)| + | A(u\rightarrow \{x_i,x_{i+1} \})| \leq 2. 
$$

\noindent\textbf {Claim 7}. If $k\in [1,m]$, then  $| A(\{x_{k-1},x_k\}\rightarrow u)| \leq 1$.

\noindent\textbf {Proof}. Suppose, to the contrary, that is $k\in [1,m]$ and  $\{ x_{k-1},x_k \}\rightarrow u$. Without loss of generality, we may assume that $A(u,x_{k+1})=\emptyset$. To be definite, assume that  $x_{k+2}:=x_1$ and $x_m:=x_{k+1}$. Then  $ux_1\notin D$ by Claim 6.

First suppose that $x_1u\in D$. It is easy to see that  $p\geq 6$ and  $ A(u\rightarrow \{x_{m-1},x_m,x_1,x_2\})=\emptyset$.\,Using this together with $od(u)\geq p/2-1$ we see that  $ A(u\rightarrow \{x_3,x_4,\ldots , x_{m-2}\})\not =\emptyset$  and  $d(u,\{x_1,x_2,\ldots ,x_{m-1}\})\geq p-3$. Note that $m\geq 6$ and show that for each  $j\in [3,m-3]$,
 $$
 |A(u\rightarrow \{x_j,x_{j+1}\}| \leq 1. \eqno {(5)}
$$

Assume that (5) is not true. Then  $u\rightarrow \{ x_j,x_{j+1}\}$ for some $j\in[3,m-3]$. We can assume that $j$ is as small as possible. Then  $A(u,x_{j-1})=\emptyset$  and  $x_{j-2}u \notin D$ by Claim 6. Hence $j\geq 4$. Since the vertex $u$ cannot be inserted into the cycle $C$, $ux_1 \notin D$  and  $x_{j-2}u \notin D$, by Lemma 2 we have
$$
d(u,\{x_1,x_2,\ldots ,x_{j-2}\})\leq j-3 \quad \hbox {and} \quad d(u,\{x_j,x_{j+1},\dots ,x_{m-1}\})\leq m-j+1.
$$
Hence $d(u)\leq p-2$, a contradiction, which proves (5).

From $A(u\rightarrow \, \{x_{m-1},x_m,x_1,x_2\})=\emptyset$ \, and (5) it follows that \, $od(u)\leq p/2-2$,\, a contradiction.

So suppose next that $x_1u\notin D$. Then $A(u,x_1)=\emptyset$ by Claim 6, $m\geq 4$ and $d(u,\{x_2,x_3,\ldots,x_{m-1}\})$ $\geq p-3$. Hence,  $ ux_2\in D$ by Lemma 2(ii). Note that $A(v,x_m)=\emptyset$ and $vx_1\notin D$. By Lemma 2(iii), it is easy to see that  $x_{m-1}v\in D$ and  $d(v,\{x_1,x_2,\ldots ,x_{m-1}\})=p-3$. 
If $x_1v\notin D$, then $A(v,x_1)=\emptyset$, and by Lemma 2, $vx_2\in D$ . Now we have  $x_{m-1}u$,\,$ vx_2\in D$ and  $A(\{u,v\}, \{x_m,x_1\})=\emptyset$, and the considered Case 1 ($h\geq 1$) holds. So we may assume that this is not the case. Then $x_1v\in D$. We also can assume that $x_{m-2}v\notin D$ (otherwise $\{ x_{m-2},x_{m-1},x_1 \}\rightarrow v$ and for the vertex $v$ the considered case $x_1u\in D$ holds). From $x_{m-2}v\notin D$ and $vx_1\notin D$, by Lemma 2(iii), it follows that $d(v,\{x_1,x_2,\ldots ,x_{m-2}\})\leq p-5$. Hence $vx_{m-1}\in D$. From $A(x_m,\{u,v\})=\emptyset$ and $d(x_m)\geq p-1$  we have $d(x_m,\{x_2,x_3,\ldots ,x_{m-1}\})\geq p-3$. Note that  $x_m$ cannot be inserted into the path $x_2x_3\ldots x_{m-1}$ (otherwise we obtain a cycle of length $m$, which does  not contain the vertices $v$ and $x_1$ and therefore, by Claim 2, the subdigraph $\langle \{v,x_{1} \} \rangle$ is strong, which contradicts the fact that $vx_1\notin D$). It follows that  $x_mx_2\in D$ by Lemma 2, and  $C_{m+1}=x_{m-2}uvx_{m-1}x_mx_2 \ldots x_{m-2}$, a contradiction. Claim 7 is proved. \fbox \\\\

Similarly to Claim 7, we can show the following:\\

\noindent\textbf {Claim 8}. If  $i\in [1,m]$, then $| A(u\rightarrow \{x_i,x_{i+1}\})| \leq 1$. \fbox \\\\

\noindent\textbf {Claim 9}. If  $k\in [1,m]$, then $| A(x_k \rightarrow u ) | +| A(u \rightarrow x_{k-1})| \leq 1$.

\noindent\textbf {Proof}. Suppose, to the contrary, that is $k\in [1,m]$ and  $x_ku$, $ux_{k-1}\in D$. To be definite, assume that  $x_k:=x_2$. From Claims 7, 8 and (1) it follows that
$$
A(u,\{x_m,x_3\})=A(x_1\rightarrow u)=A(u\rightarrow x_2)=\emptyset.
$$
From this it is easy to see that $m\geq 5$ \, and $d(u,\{x_4,x_5,\ldots ,x_{m-1})\geq p-5$. Since the vertex $u$ cannot be inserted into the path $x_4x_5\ldots x_{m-1}$, by Lemma 2  we have $ux_4$ and $x_{m-1}u\in D$. Hence,  $| A(v,x_3)| =| A(v,x_m)| =2$  by Claim 5. Therefore $A(v,\{x_1,x_2\})=\emptyset$ by Claims 7 and 8. Since $C_{m+1}\not\subset D$ it is not difficult to see that $x_mx_2$,\,$x_2x_4$,\,$x_1x_3$ \,and $x_3x_2 \,\notin D$ (if $x_3x_2\in D$, then $C_{m+1}=x_mvx_3x_2ux_4\ldots x_m$). So we have \,$d(x_2,\{u,v,x_1,x_3\})\leq 4$ and $d(x_2,\{x_4,x_5,\ldots ,x_m\}\geq p-5$. Therefore, since $ x_mx_2\notin D$ and $x_2x_4\notin D$, applying Lemma 2(iii), we can insert $x_2$  into the path $x_4x_5\ldots x_m$ (i.e.,   $x_ix_2$,  $x_2x_{i+1} \in D$ for some $i\in [4,m-1]$) and obtain a cycle $x_mvux_4\ldots x_ix_2x_{i+1}\ldots x_m$ \, of length $m$, which  does not contain the vertices $x_1$ and $x_3$. By Claim 2, the subdigraph $\langle \{x_1,x_3 \} \rangle$ is strong. Hence $x_1x_3\in D$,  which contradicts the fact that $x_1x_3\notin D$. This completes the proof of Claim 9. \fbox \\\\

We now divide Case 2 ($h=0$) into two subcases.\\

\noindent\textbf {Subcase 2.1}. \, $ p=2n+1$.

From  (1) and Claims 7- 9 it follows that the vertex $u$ (respectively, $v$) is  adjacent  with at most one vertex of two consecutive vertices of the cycle $C$ and  $ O(u,C)=I(u,C)$ and $O(v,C)=I(v,C)$. Hence, without loss of generality, we may assume that 
$$
A(u,\{x_2,x_3\})=\emptyset \quad \hbox {and} \quad O(u)=I(u)=\{x_1,x_4,x_6,\ldots , x_{p-3},v \}.  \eqno {(6)}
$$

If $m=3$, then  Claim 3 implies that $| A(v,x_2)| =2$ or $|A(v,x_3)|=2$ and $C_4\subset D$, a contradiction.

Assume that $m=2n-1\geq 5$. Since $x_{m-1}u\in D$, by Claim 5 we have  $| A(x_m,v)| =2 $. Therefore, by an argument similar to (6), we get that ether  $| A(v,x_2)| =2$ or $|A(v,x_3)| =2$. From this and (6) it is easy to see that if $vx_3\in D$,  then  $C_{m+1}=x_1uvx_3x_4\ldots x_mx_1$ and if $x_2v\in D$,  then  $C_{m+1}=x_2vux_4x_5\ldots x_mx_1x_2$, a contradiction.\\

 \noindent\textbf {Subcase 2.2.}  $ p=2n $.

From $d(u)\geq 2n-1$ it follows that ether $od(u)\geq n$ or $id(u)\geq n$. Without loss of generality, we may assume that $od(u)\geq n$ (otherwise we will consider the digraph $\overleftarrow {D}$). Now from Claims 7-9 it follows that
$$
u\rightarrow \{x_1,x_3,\ldots ,x_{2n-3}\}, \quad A(u,\{x_2,x_4,\ldots ,x_{2n-2}\})=\emptyset, \eqno {(7)}
$$
$$
I(u)\subseteq \{v,x_1,x_3,\ldots ,x_{2n-3}\}.  \eqno {(8)}
$$
Since $id(u)\geq n-1$, without loss of generality,  we may assume that $\{x_1,x_3,\ldots ,x_{2n-5}\}\rightarrow u$. Hence , by (7) and Claim 5, it follows that for each $i\in [1,n-2]$,
$$
| A(u,x_{2i-1})| =| A(v,x_{2i})|=2.      \eqno {(9)}
$$
Then by Claims 7-9 we have $A(v,\{x_1,x_3,\ldots ,x_{2n-3}\})=\emptyset$. Therefore  $A(v,x_{2n-2})\not=\emptyset$ since $d(v)\geq 2n-1$, i.e. $vx_{2n-2}\in D$ or $v_{2n-2}v\in D$. If $vx_{2n-2}\in D$, then  $x_{2n-3}u\in D$ and $x_{2n-2}v\in D$ by Claim 5, (8) and (9). So, in any case we have that $x_{2n-2}v\in D$. Then $id(v)\geq n$.\\ 

We will now show that
$$
A(\langle \{ x_1,x_3,\ldots ,x_{2n-3}\} \rangle )=\emptyset . \eqno {(10)}  
$$

\noindent\textbf {Proof of (10)}. Assume that  (10) is not true. Then $x_ix_j\in D$ for some distinct vertices $x_i, x_j \in \{ x_1,x_3,\ldots ,$ $x_{2n-3} \}$. Assume   that  $| \{x_{i+1},x_{i+2},\ldots ,x_{j-1}\}| =1$. Then from (9) and  $x_{2n-2}v\in D$ we have: a) if $j=2n-3$, then  $i=2n-5$ and   $C_{m+1}=x_{2n-5}x_{2n-3}x_{2n-2}vux_1x_2\ldots x_{2n-5}$; b) if $j\not= 2n-3$, then  $C_{m+1}=x_ix_juvx_{j+1}\ldots x_{i-1}x_i$. Now assume that  $| \{x_{i+1},x_{i+2},\ldots ,x_{j-1}\}| \geq 2$. Then $n\geq 6$. Using (9) we can see that: c) if $i\not= 2n-3$ and $j\not= 1$, then $C_{m+1}=x_ix_jx_{j+1}\ldots x_{i-1}vx_{i+1}\ldots x_{j-2}ux_i$; d) if $i=2n-3$ or $j=1$, then $C_{m+1}=x_ix_jx_{j+1}\ldots x_{i-2}ux_{i+2}\ldots x_{j-1}vx_{i-1}x_i$.  Hence in each case we have a $C_{m+1}\subset D$, which is a contradiction and  (10) is proved. \fbox \\\\

 Using an analogous argument for $\overleftarrow {D}$, similarly to (10), we can show that
$
A(\langle \{ x_2,x_4,\ldots , x_{2n-2}\} \rangle )=\emptyset .
$
Therefore 
$$
A(\langle \{v, x_1,x_3,\ldots , x_{2n-3}\} \rangle )=A(\langle \{u, x_2,x_4,\ldots , x_{2n-2} \} \rangle )=\emptyset \quad \hbox {and} \quad D\subseteq K_{n,n}^*.
$$
This contradicts the our supposition. The discussion of Case 2 is completed and  Theorem 1 is proved. \fbox \\\\

\noindent\textbf {3. A sufficient condition for a digraph to be hamiltonian}\\
 
\noindent\textbf{Theorem 2}. Let $D$ be a digraph on $2n\geq 6$ vertices with minimum degree  at least $2n-1$ and with minimum semi-degree  at least $n-1$. Then $D$ is hamiltonian unless  
$$
D\in H(n,n)\cup H(n,n-1,1) \cup \{H(2n),H'(2n),D_6,D_6',\overleftarrow D_6, \overleftarrow {D_6'}\}.
$$

 \noindent\textbf{Proof}. By Lemma 4(i), the result is easily verified if $D$ is not strong. Now assume that $D$ is strong. The proof is by contradiction. Suppose that Theorem 2 is false, in particular, $D$ is not hamiltonian. Then it is not difficult to see that $D\not\subseteq  K_{n,n}^*$. By Theorem 1, $D$ has a cycle of length $2n-1$. Let $C:= C_{2n-1}:=x_1x_2\ldots x_{2n-1}x_1$  be an arbitrary cycle of length $2n-1$ in $D$ and let the vertex $x$ is not containing  this cycle $C$. Since $C$ is a longest cycle, using Lemmas 1 and 2, we obtain the following claim:\\

\noindent\textbf {Claim 1.} (i) $d(x)=2n-1$ and there is a vertex $x_l$, $l\in [1,2n-1]$  such that $A(x,x_l)=\emptyset$.

 (ii) If $x_ix\notin D$, then $xx_{i+1}\in D$ and if $xx_i\notin D$, then $x_{i-1}x\in D$, where $i\in [1,2n-1]$.

(iii) If $A(x,x_i)=\emptyset$, then $x_{i-1}x,\,xx_{i+1}\in D$ and $d(x_i)=2n-1$. \fbox \\\\

 By Claim 1(i), without loss of generality, we may assume that $A(x,x_{2n-1})=\emptyset$. For convenience, let $p:=2n-2$ and $y:=x_{2n-1}$. We have $ yx_1,$ $x_py\in D$ and  $x_px$, $xx_1\in D$ by Claim 1(iii), and $d(y)=2n-1$ by Claim 1(i).

Let $\{u,v\}:=\{x,y\}$ and for each $z\in \{x,y\}$ let  
 $$
O^-(z):=\{x_i/zx_{i+1}\in D,\, i\in [1,p-2] \},\quad I^+(z):=\{x_i/ x_{i-1}z\in D,\, i\in [2,p-1] \}.
$$
 We first prove the following Claims 2-11.\\

\noindent\textbf{Claim 2}. If $x_{p-1}u \in D$, then $A(O^-(v)\rightarrow x_p)=\emptyset$.

 \noindent\textbf{Proof}. Assume, to the contrary, that  $x_{p-1}u \in D$ and $x_ix_p \in D$, where $x_i\in O^-(v)$. Then by the definition of $O^-(v)$, $vx_{i+1}\in D$ and $x_1x_2\ldots x_ix_pvx_{i+1}$ $\ldots x_{p-1}ux_1$ is a hamiltonian cycle, a contradiction. \fbox \\\\
  
\noindent\textbf{Claim 3}. If $x_{p-1}u,\,vx_p\in D$, then $A(x_p \rightarrow I^+(v))=\emptyset$.

 \noindent\textbf{Proof}.  Assume, to the contrary, that  $x_{p-1}u$, $vx_p$ and $x_px_i\in D$, where $x_i\in I^+(v)$. Then by the definition of $I^+(v)$, $x_{i-1}v\in D$ and  $x_1x_2\ldots x_{i-1}vx_p$ $x_i\ldots x_{p-1}ux_1$ is a hamiltonian cycle, a contradiction. \fbox \\\\

\noindent\textbf{Claim 4}. If $x_{p-1} \rightarrow \{x,y\}$, then $od(x)=od(y)=n-1$, $id(x)=id(y)=n$  and $O(x)=O(y)$.

\noindent\textbf{Proof}. Let $ x_{p-1}\rightarrow \{x,y\}$. Since $C$ is longest cycle of $D$, we have $A(\{x,y\}\rightarrow x_p)=\emptyset$. By Claim 2,  $A(O^-(x)\cup O^-(y)\rightarrow x_p)=\emptyset$. Hence,  $| O^-(x)\cup I^-(y)\cup \{x,y\}| \leq n$ by Lemma 4(ii). Therefore, since  $|O^-(u)|=od(u)-1$, we deduce that $od(x)=od(y)=n-1$ and  $O(x)=O(y)$. Together with  $d(x)=d(y)=2n-1$ this implies that  $id(x)=id(y)=n$. Claim 4 is proved. \fbox \\\\

Similarly to Claim 4, we can show the following claim:\\  

\noindent\textbf{Claim 5}. If $\{x,y\} \rightarrow x_2$, then  $id(x)=id(y)=n-1$, $od(x)=od(y)=n$ and $I(x)=I(y)$. \fbox \\\\

\noindent\textbf{Claim 6}.  $| A(x_{p-1}\rightarrow \{x,y\})| \leq 1$.

\noindent\textbf{Proof}. Assume, on the contrary, that $x_{p-1}\rightarrow \{x,y\}$. Then  $id(x)=id(y)=n$,  $ od(x)=od(y)=n-1$ and  $O(x)=O(y)$ by Claim 4. Hence, $xx_2\in D$ if and only if $yx_2\in D$. Therefore, $A(\{x,y\}\rightarrow x_2)=\emptyset$ by Claim 5. Hence,  $x_1\rightarrow \{x,y\}$ by Claim 1(ii). Together with $od(x)\geq 2$ this implies that there exists an  $k\in [2,p-2]$ such that $x_{k-1}x$,  $xx_{k+1}\in D$  and  $A(x,x_k)=\emptyset$. Applying Claim 1(iii)  we find that $d(x_k)=2n-1$. From $O(x)=O(y)$ it is not difficult to see that $A(x_k,y)=\emptyset$, $yx_{k+1}$, $x_{k-1}y\in D$. Then by Lemma 2, since $x_k$ cannot be inserted into the path $x_1x_2\ldots x_{k-1}$ and into the path $x_{k+1}x_{k+2}\ldots x_p$,  we have
 $$
 d(x_k,\{x_1,x_2,\ldots , x_{k-1}\})\leq k \quad \hbox { and} \quad d(x_k,\{x_{k+1},x_{k+2},\ldots ,x_p\})\leq p-k+1.
$$
Using this,  $d(x_{k})=p+1$,   $A(x_k,\{x,y\})=\emptyset$ and  Lemma 2(ii), we obtain     
$$
 x_kx_1, x_px_k \in D \quad \hbox {and} \quad  d(x_k,\{x_1,x_2,\ldots ,x_{k-1}\})=k, \, d(x_k,\{x_{k+1},x_{k+2},\ldots , x_p\})=p-k+1.  \eqno {(11)}
$$

We now show that 
$$
A(\{x,y\}\rightarrow x_{k+2})=\emptyset. \eqno {(12)}
$$

\noindent\textbf{Proof of (12)}. Suppose that (12) is false. From $O(x)=O(y)$ it is clear that  $\{x,y\}\rightarrow x_{k+2}$. Since $x_kx_1\in D$ we see that $x_{k+1}x_k\notin D$ (otherwise  $x_{k+1}x_k\in D$ and  $x_1x_2\ldots x_{k-1}xx_{k+2}\ldots x_pyx_{k+1}x_kx_1$  is a hamiltonian cycle, a contradiction). Note that  $x_1x_2\ldots x_{k-1}xx_{k+1}\ldots x_pyx_1$ is a cycle of length $2n-1$ and the vertex $x_k$ cannot be inserted into this cycle. Then, since $x_{k+1}x_k\notin D$ and $d(x_k)=2n-1$, using Claim 1(ii) we get  $x_kx_{k+2}\in D$. From this it follows that $\{x,x_k\}\rightarrow \{x_{k+1},x_{k+2}\}$  and  $x_{k-1}\rightarrow \{x,x_k\}$ for the path  $x_{k+1}x_{k+2}\ldots x_pyx_1x_2\ldots x_{k-1}$. Therefore $id(x)=n-1$ by Claim 5, which  contradicts the fact that $id(x)=n$. This proves (12). \fbox \\

From $x_1\rightarrow \{x,y\}$ and (12) it follows that  $|A(x\rightarrow \{x_i,x_{i+1}\})| \leq 1$ for each $i\in [1,p-2]$. Therefore, since $od(y)=od(x)=n-1$  and $O(x)=O(y)$,  it is not difficult to see that
$$
\{x,y\}\rightarrow \{x_1,x_3,\ldots ,x_{p-1}\}\rightarrow \{x,y\}, \eqno {(13)}
$$
$$
A(\{x,y\},\{x_2,x_4,\ldots ,x_{p-2}\})=\emptyset. \eqno {(14)}
$$
Together with Claim 2  this implies that
$$
A(\{x_2,x_4,\ldots ,x_{p-2}\}\rightarrow x_p)=\emptyset. \eqno {(15)}
$$

It is not difficult to show that
$$
A(\langle \{x_2,x_4,\ldots ,x_{p-2}\} \rangle )=\emptyset. \eqno {(16)}
$$

Indeed, if (16) is false, then $x_ix_j\in D$ for some distinct vertices $x_i, x_j \in \{x_2,x_4,\ldots ,x_{p-2}\}$. It is easy to see that if $i<j$, then $C_{2n}=x_1x_2\ldots x_ix_j\ldots x_pyx_{i+1}\ldots x_{j-1}xx_1$ and if $i>j$,  then $x_jx_1\in D$ by (11), and  $C_{2n}=x_1x_2\ldots x_{j-1}xx_{i+1}\ldots x_pyx_{j+1}\ldots x_ix_jx_1$, a contradiction.

From (14) and (16) it follows that $A(\langle \{x,y,x_2,x_4,\ldots ,x_{2n-4}\} \rangle )=\emptyset.$
 By (15), now it is not difficult to see that  $D\in H(n,n-1,1)$, where $a:=x_p$, $ A:=\{x,y,x_2,x_4,\ldots ,x_{2n-4}\}$ and  $B:=\{x_1,x_3,\ldots ,$ $x_{2n-3}\}$. This contradicts to our supposition that  $D\notin H(n,n-1,1)$. The proof of Claim 6 is completed. \fbox \\\\

Similarly to Claim 6, we can show the following claim:\\

\noindent\textbf{Claim 7}. $| A(\{x,y\}\rightarrow x_2)| \leq 1$. \fbox \\\\

\noindent\textbf{Claim 8}. There is a vertex $x_i$, $i\in [2,p-1]$, such that $A(x,x_i)=\emptyset$ (i.e., if $C_{2n-1}$, $p=2n-2$, is an arbitrary cycle of $D$ and the vertex $x\notin C_{2n-1}$, then $x$ is not adjacent with at least two vertices).

\noindent\textbf{Proof}. Suppose, on the contrary, that the vertex $x$ is adjacent with each vertex $x_i$,  $i\in [2,p-1]$. Since $n\geq 3$, $d(x)=p+1$ and $D$ is not hamiltonian, there is an  $l\in [2,p-1]$  such that 
$$
O(x)=\{x_1,x_2,\ldots ,x_l\} \quad \hbox {and } \quad I(x)=\{x_l,x_{l+1},\ldots ,x_p\}. \eqno {(17)}
$$
 Since  $od(x)$  and  $id(x)\geq n-1$, we see that $l=n-1$ or $ l=n$. Hence  $x_{p-1}y\notin D$ and $yx_2\notin D$ by Claims 6 and 7. Now  $x_1y$ and $yx_p\in D$ by Claim 1(ii). Therefore, since $y$ cannot be inserted into the path $x_1x_2\ldots x_p$, there is a vertex $x_k$, $k\in [2,p-1]$, such that $A(y,x_k)=\emptyset$. Using Claim 1(iii),  we get 
$$
x_{k-1}y, \, yx_{k+1}\in D \quad \hbox {and} \quad d(x_k)=p+1. \eqno {(18)}
$$ 
Choose $k$ is as large as possible. It follows that $y\rightarrow \{x_{k+1}, x_{k+2},\ldots ,x_p\}$ and $k\geq l-1$. We can assume that $k\geq l$ (if $k=l-1$, then in digraph $\overleftarrow D$ we will have the case $k\geq l+1$).

Suppose first that $k\geq l+1$. If $x_ix_k\in D$, where $i\in [1,l-1]$, then by (17) and (18), $x_1x_2\ldots x_ix_k\ldots x_p$ $xx_{i+1}\ldots x_{k-1}yx_1$ is a hamiltonian cycle, a contradiction. So we may assume that $A(\{x_1,x_2,\ldots ,x_{l-1}\}\rightarrow x_k)=\emptyset$. Using this together with   $A(\{x,y\}\rightarrow x_k)=\emptyset$, $l\geq n-1$ and $id(x_k)\geq n-1$, we obtain $x_px_k\in D$. Therefore $x_1x_2\ldots x_{k-1}yx_{k+1}\ldots x_px_kxx_1$ is a hamiltonian cycle, a contradiction.

Now  suppose that $k=l$. Assume, without loss of generality, that $A(x_i,y)\not=\emptyset$ for each $i\in [2,l-1]$ (otherwise in  $\overleftarrow D$ we will have the considered case $k\geq l+1$). Then from $x_1y\in D$ it follows that
$$
 \{x_1,x_2,\ldots ,x_{l-1}\}\rightarrow y.  \eqno {(19)}
$$
We also can assume that $l=n$ (if $l=n-1$, then in  $\overleftarrow D$ we will have the case $l=n$). So, we have $k=l=$ $n$. It is not  difficult to see that 
$$
A(x_1,x_n)=  
A(\{x_1,x_2,\ldots ,x_{n-2},x_p\}\rightarrow x_n)=A(x_n\rightarrow \{x_{n+2},x_{n+3},\ldots ,x_{p-1}\})=\emptyset. \eqno {(20)}
$$
Indeed, if it is not true, then  by (17) and (19) we have

   if $x_ix_n\in D$ and $i\in [1,n-2]$, then $C_{2n}=x_1x_2 \ldots  x_ix_n \ldots x_pxx_{i+1} \ldots x_{n-1}yx_1$; 

   if $x_nx_i\in D$ and $i\in [n+2,p-1]$, then \ $C_{2n}=x_1x_2 \ldots x_nx_ix_{i+1}\ldots x_pyx_{n+1}\ldots x_{i-1}xx_1$; 

   if $x_px_n \in D$, then  $C_{2n}=x_1x_2 \ldots x_{n-1}yx_{n+1} \ldots x_px_nxx_1$;

   if $x_nx_1 \in D$, then  $C_{2n}=x_1x_2 \ldots x_{n-1}yx_{n+1} \ldots x_pxx_nx_1$. In each case we have a  hamiltonian cycle, a contradiction, and (20) holds.   

 Therefore from $d(x_n)=2n-1$  and (20), since $x_n$ cannot be inserted into the paths $x_1x_2\ldots x_{n-1}$  and $x_{n+1}x_{n+2}\ldots x_p$, it follows that (by Lemma 2)
$$
\{x_{n+1},x_{n+2},\ldots ,x_{p-1}\}\rightarrow x_n\rightarrow \{x_2,x_3,\ldots ,x_{n-1}\}. \eqno {(21)}
$$
If $x_ix_1\in D$ for some $i\in [2,p-1]\setminus\{n\}$, then by (17), (18), (19) and (21) we have  if $i\in [2,n-1]$, then $C_{2n}=x_1x_2\ldots x_{i-1}yx_{n+1}\ldots x_pxx_{i+1}\ldots x_nx_ix_1$ and if  $i\in [n+1,p-1]$, then $C_{2n}=x_1x_2\ldots x_{n-1}yx_{i+1}\ldots x_pxx_n$ $\ldots x_ix_1$, a contradiction. So, we may assume that
$$
A(\{x_2,x_3,\ldots ,x_{p-1}\}\rightarrow x_1)=\emptyset. \eqno {(22)}
$$ 
Hence, by Lemma 4(ii), $2n-4\leq n$, i.e. $n\leq 4$. Let $n=4$. Then by (22,) $id(x_1)\leq 3$. On the other hand, from $A(x_1\rightarrow \{x,x_3,x_4,x_5\})=\emptyset$ it follows (if $x_1x_5\in D$, then $x_4x_2\in D$ by (21), and $C_8=x_1x_5x_6xx_3x_4x_2yx_1$) that $od(x_1)\leq 3$. So $d(x_1)\leq 6$, a contradiction. Let now $n=3$. From (21) we see  that $x_3x_2\in D$. Hence it is easy to see that $x_4x_3\notin D$ by (20), $x_4x_2 \notin D$ and
$$
A(x_1 \rightarrow \{x,x_3,x_4 \})=A(\{x_2,x_3\}\rightarrow  x_1)=\emptyset.
$$
Therefore $x_4x_1\in D$. Now, it is not difficult to check that $D$ is isomorphic to one of the digraphs $D_6$, $D_6'$, a contradiction. This completes the proof of Claim 8. \fbox \\\\

Similarly to Claim 8, we can show the following claim:\\

\noindent\textbf{Claim 9}. There is a vertex $x_i$, $i\in [2,p-1]$, such that $A(y,x_i)=\emptyset$.  \fbox \\\\

\noindent\textbf{Claim 10}.  $x_{p-1}y\notin D$.

\noindent\textbf{Proof of claim 10}. Suppose, on the contrary, that $x_{p-1}y\in D$. By Claim 6 we have $x_{p-1}x\notin D$. Therefore  $xx_p\in D$ by Claim 1(ii) . By Claim 8 there is a vertex $x_l$, $l\in [2,p-1]$, such that  $A(x,x_l)=\emptyset$. Using Claim 1(iii), we obtain 
$$
x_{l-1}x,\, xx_{l+1}\in D \quad \hbox {and} \quad d(x_l)=2n-1=p+1. \eqno {(23)}
$$
 For the vertex $x_l$ we first will prove the following statements \textbf{(a)-(i)}. \\

\noindent\textbf{(a).} $x_px_l\notin D$.

\noindent\textbf{Proof}. Indeed, if \textbf{(a)} is not true, then $x_px_l\in D$ and  $C_{2n}=x_1x_2\ldots x_{l-1}xx_px_l\ldots x_{p-1}yx_1$ by (23),  a contradiction. \fbox \\\\

\noindent\textbf{(b)}. If $l\leq p-2$, then $A(x_l,x_p)=\emptyset$.

\noindent\textbf{Proof}. From $l\leq p-2$ it follows that $x_l\in O^-(x)$. Hence $x_lx_p\notin D$  by Claim 2. Therefore by  statement \textbf{(a)}, $A(x_l,x_p)=\emptyset$. \fbox \\\\

\noindent\textbf{(c)}. If $l\leq p-2$, then $x_{p-1}x_l$ and  $x_ly\in D$.

\noindent\textbf{Proof}. Note that  $A(x_l,x_p)=\emptyset$ by statement \textbf{(b)}, and  the cycle $x_1x_2\ldots x_{l-1}xx_{l+1}\ldots x_pyx_1$ has length $2n-1$. Therefore   $x_{p-1}x_l$ and  $x_ly\in D$ by Claim 1(iii). \fbox \\\\

\noindent\textbf{(d)}. If $l\leq p-2$, then $A(y\rightarrow \{x_l,x_{l+1},\ldots ,x_p\})\not= \emptyset$.

\noindent\textbf{Proof}. Suppose, on the contrary, that $A(y\rightarrow \{x_l,x_{l+1},\ldots ,x_p\})= \emptyset$. It follows that $O^-(y)\subseteq \{x_1,x_2,\ldots ,x_{l-2}\}$. If $x_i\in O^-(y)$ and $x_ix_l\in D$, then $C_{2n}=x_1x_2\ldots x_ix_l\ldots x_pyx_{i+1}\ldots x_{l-1}xx_1$ is a hamiltonian cycle in $D$, a contradiction. So we can assume that $A(O^-(y)\rightarrow x_l)=\emptyset$. Together with  $A(\{x,y\}\rightarrow x_l)=\emptyset$ and $|O^-(y)| \geq n-2$ this implies that $x_px_l\in D$. But this contradicts \textbf{ (a)}, and  hence \textbf{(d)} is proved. \fbox \\\\

\noindent\textbf{(e)}. If $l\leq p-2$, then $x_px_{l+1}\notin D$ and $x_{l-1}x_p\notin D$.

\noindent\textbf{Proof}. Recall that $xx_p,\, x_px\in D$, and $x_{p-1}x_l,\, x_ly\in D$ by \textbf{(c)}. Then by (23) we have, if $x_px_{l+1}\in D$, then  $C_{2n}=x_1x_2\ldots x_{l-1}xx_px_{l+1}\ldots x_{p-1}x_lyx_1$ and if $x_{l-1}x_p\in D$, then $C_{2n}=x_1x_2\ldots x_{l-1}x_pxx_{l+1}\ldots $ $ x_{p-1}x_l$ $yx_1$. Therefore $D$ is hamiltonian, a contradiction. \fbox \\\\

\noindent\textbf{(f)}. If $l\leq p-2$ and $xx_{l+2}\in D$, then $x_lx_{l+2}\notin D$ and $x_{l+1}x_l\in D$.

\noindent\textbf{Proof}. Indeed, if $x_lx_{l+2}\in D$, then for the path $x_{l+1}x_{l+2}\ldots x_pyx_1x_2\ldots x_{l-1}$ we have $\{x,x_l\}\rightarrow \{x_{l+1},x_{l+2}\}$ and $x_{l-1}\rightarrow \{x,x_l\}$, which contradicts  Claim 7. So $x_lx_{l+2}\notin D$. Now from Claim 1(ii) it follows that $x_{l+1}x_l\in D$. \fbox \\\\

\noindent\textbf{(g)}. If $l\geq 3$ and $x_{l-2}x\in D$, then $x_{l-2}x_l\notin D$ and $x_lx_{l-1}\in D$.

\noindent\textbf{Proof}. Indeed, if $x_{l-2}x_l\in D$, then for the path $x_{l+1}x_{l+2}\ldots x_pyx_1\ldots x_{l-2}x_{l-1}$ we have $\{x_{l-2},x_{l-1}\}$ $\rightarrow \{x,x_l\}$ and $\{x,x_l\}\rightarrow x_{l+1}$, which contradicts Claim 6. So  $x_{l-2}x_l\notin D$. From this and Claim 1(ii) it follows that $x_lx_{l-1}\in D$. Statement \textbf{(g)} is proved. \fbox  \\\\

 \noindent\textbf{(h)}. If $l\leq p-2$, then $x_ix_p\in D$ if and only if $ x_i\notin \{x_{l-1}\}\cup O^-(x)$;  and  $x_px_i\in D$ if and only if $x_i\notin \{x_{l+1}\}\cup I^+(x)$.

\noindent\textbf{Proof}.  By Claims 2, 3 and statement \textbf{(e)} we have
$$
A(O^-(x)\cup \{x_{l-1},y\}\rightarrow x_p)=A(x_p\rightarrow \{x_{l+1}\}\cup I^+(x))=\emptyset.  \eqno (24)
$$
From $x_{p-1}x\notin D$ and $xx_p\in D$, we get that 
$$
| I^+(x)| =id(x)-1 \quad  \hbox {and} \quad | O^-(x)| =od(x)-2. 
$$ 
Therefore   $id(x_p)\leq  2n-1-od(x)$  and   $ od(x_p) \leq  2n-1-id(x)$  by (24).  Hence $id(x_p)=  2n-1-od(x)$  and   $ od(x_p) = 2n-1-id(x)$ (otherwise   $ d(x)+d(x_p)<4n-2$, which is a contradiction). Now from this it is not difficult to see that  statement \textbf{(h)} is true. \fbox \\\\

Recall that the proof of statement \textbf{(h)} implies the following statement:\\

\noindent\textbf{(i)}. The vertex $x$ is not adjacent with at most one vertex of the path $x_1x_2\ldots x_{p-2}$, in particular,  the vertex $x$ is not adjacent with at most  3 vertices (i.e., if $C_{2n-1}$ is an arbitrary cycle of $D$ and the vertex $x\notin C_{2n-1}$, then $x$ is not adjacent with at most  tree vertices).  \fbox \\\\

By Claim 8 there is a vertex $x_k$, \, $k\in [2,p-1]$, such that $A(x,x_k)=\emptyset$. Without loss of generality, assume that $k$ is as large as possible. From the maximality of $k$ and Claim 1(iii) it is easy to see that   
$$
x_{k-1}x\in  D, \, d(x_k)=p+1,\,x\rightarrow \{x_{k+1},x_{k+2},\ldots ,x_p\}, \,   A(\{x_{k+1},x_{k+2},\ldots ,x_{p-1}\}\rightarrow x)=\emptyset . \eqno {(25)}
$$

 We now consider  two cases. 

\noindent\textbf{Case 1}. $k\leq p-2$.

Then by statement \textbf{(c)} we have
$$
x_{p-1}x_k \in D \quad \hbox {and} \quad x_ky \in D. \eqno {(26)} 
$$

From statement \textbf{(i)} and (25) it follows that if $i\in [1,p]$ and $i\not= k$, then
$$
A(x,x_i)\not= \emptyset. \eqno {(27)}
$$

It is easy to see that $n\geq 4$. Indeed,  if $n=3$, then $k=2$ and by (26) the vertex $y$ is not adjacent only with  one  vertex of  the cycle $x_1x_2\ldots x_pxx_1$, which contradicts Claim 9.

Suppose first that $k\leq p-3$. Then $x_1\notin I^+(x)$, and  $x_{k+2}\notin I^+(x)$ by (25). Together with   statement \textbf{(h)} this implies that 
$$
x_p \rightarrow \{ x_1,x_{k+2} \}.  \eqno {(28)}
$$

If $x_{k+1}y\in D$, then using  (25), (26) and (28), we obtain $C_{2n}=x_1x_2\ldots x_{k-1}xx_px_{k+2}\ldots x_{p-1}x_kx_{k+1}y$ $x_1$, a contradiction. So, we may assume that $x_{k+1}y\notin D$. Since  $x_ky\in D$ by (26), we see that  $A(y,x_{k+1})=\emptyset$. Therefore  $yx_{k+2}\in D$ by Claim 1(iii). Recall that $x_{k+1}x_k \in D$ by statement \textbf{(f)}, and hence by (25) and (28) we have a hamiltonian cycle $x_1x_2\ldots x_{k-1}xx_{k+1}x_kyx_{k+2}\ldots x_px_1$, a contradiction.

Suppose next that $k=p-2$. Then by  $x_{p-1}y\in D$ and statements \textbf{(d)},  \textbf{(c)},  $yx_{p-2}\in D$.
If $xx_2\in D$, then $x_1x\notin D$, $x_2\notin I^+(x)$ and  $x_px_2\in D$ by statement \textbf{(h)}. Since $xx_2\in D$, by Claim 7 we have $yx_2\notin D$. Therefore $x_1y\in D$  by Claim 1(ii), and we get a hamiltonian cycle $x_1yx_{p-2}x_{p-1}x_px_2\ldots x_{p-3}xx_1$, a contradiction. So we may assume that $xx_2\notin D$. From this and (27) it follows that
  
$$\{x_1,x_2,\ldots ,x_{p-3}\}\rightarrow x \quad \hbox {and} \quad  A(x\rightarrow \{x_2,x_3,\ldots ,x_{p-3}\})=\emptyset.$$

Therefore  $n=4$ (i.e., $p=6$). Then $x_4x_3\in D$ by statement \textbf{(g)}. We can assume that $yx_2\notin D$ (otherwise  $yx_2\in D$ and for $\overleftarrow D$ we will have the considered case $k\leq p-3$).  From Claim 1(ii) and $od(y)\geq 3$ it is easy to see that $x_1y,\, yx_3\in D$ and $A(y,x_2)=\emptyset$. Since $x_1\notin I^+(x)$, we have $x_6x_1\in D$ and $x_6x_2\notin D$ by statement \textbf{(h)}. Then $x_5x_2\notin D$ (otherwise $x_5x_2\in D$ and $C_{2n}=x_1yx_4x_5x_2x_3xx_6x_1$). Now we have $A(\{x,y,x_5,x_6\}\rightarrow x_2)=\emptyset$. Hence $x_4x_2\in D$ and $C_{2n}=x_1yx_3x_4x_2xx_5x_6x_1$, a contradiction.\\

\noindent\textbf{Case 2}. $k=p-1$.

Suppose first that $yx_{p-1}\notin D$. Then it is not difficult to see that $A(O^-(y)\rightarrow x_{p-1})=\emptyset$ (otherwise if $x_i\in O^-(y)$  and $x_ix_{p-1}\in D$, then  $C_{2n}=x_1x_2\ldots x_ix_{p-1}x_pyx_{i+1}\ldots x_{p-2}xx_1$). This together with $|O^-(y)| =od(y)-1$ and  $A(\{x,y,x_p\}\rightarrow x_{p-1})=\emptyset$ implies that $id(x_{p-1})\leq n-2$, a contradiction.  

Suppose next that $yx_{p-1}\in D$. We assume that $n\geq 5$ (It is tedious, but not difficult to prove the theorem in this case for $n=3$ and $4$. We leave its proof to the reader).\\

\noindent\textbf{Subcase 2.1}. $xx_2\in D$.

Then  $x_1x\notin D$. Using Claims 7 and 1(ii), we obtain $yx_2 \notin D$ and  $x_1y\in D$. We may assume that $A(x_2,y)=\emptyset$ (otherwise for the vertex $y$ in digraph $\overleftarrow D$ we have the considered Case 1 ($k\leq p-2)$). Then by Claim 1(iii), $yx_3\in D$. Similarly to $yx_{p-1}\in D$, we also may assume that $x_2x\in D$. From $n\geq 5$ it follows that $A(x,x_s)=\emptyset$ for some $s\in [3,p-3]$  (otherwise $O(x)=\{ x_1,x_2,x_p \}$, i.e., $od(x)\leq 3$, a contradiction). Since $x_2\notin \{x_{s+1}\}\cup I^+(x)$, using statement \noindent\textbf{(h)}, we see that $x_px_2\in D$ and $x_1yx_3\ldots x_px_2xx_1$ is a hamiltonian cycle, a contradiction.\\

\noindent\textbf{Subcase 2.2}. $xx_2\notin D$.

Then  $x_1x\in D$ by Claim 1(ii). By statement \textbf{(i)}, the vertex $x$ is not adjacent with at most one vertex of $\{x_1,x_2, \ldots, x_{p-3}\}$. From this and $n\geq 5$ it follows that $A(x,x_s)=\emptyset$ exactly for one $s\in [2,p-4]$ (otherwise $A(x,x_i)\not =\emptyset$ for each $i\in [2,p-4]$ and by $xx_2\in D$, $A(x\rightarrow \{x_2,x_3, \ldots, x_{p-3}\})=\emptyset$, i.e. $O(x)\subseteq \{ x_1,x_{p-2},x_p \}$ and $od(x)\leq 3$, which contradicts that $n\geq 5$). 

Let $s=2$ (i.e., $A(x,x_2)=\emptyset$). Note that  $x_1x$, $xx_3 \in D$ by Claim 1(iii). From  statement \textbf{(c)} it follows that $x_{p-1}x_2$ and $x_2y\in D$. Since $x_1\notin I^+(x)$, by statement \textbf{(h)} we have $x_px_1\in D$. If $yx_2\notin D$, then   $x_1y\in D$ by Claim 1(ii), and $x_1yx_{p-1}x_2\ldots x_{p-2}xx_px_1$ is a hamiltonian cycle, a contradiction. So we may assume that $yx_2\in D$. Also we may assume that $id(x)=n$ (for otherwise we will consider the digraph $\overleftarrow D$). Then, since $n\geq 4$, we see that  $\{x_{p-3},x_{p-2}\}\rightarrow x$ and by \textbf{(g)}, $x_{p-1}x_{p-2}\in D$. Then $x_{p-3}\notin O^-(x)\cup \{x_1\}$ and by \textbf{(h)} ($l=2$) we see that $x_{p-3}x_p\in D$. Thus $x_1x_2\ldots x_{p-3}x_pyx_{p-1}x_{p-2}xx_1$ is a hamiltonian cycle, a contradiction.

  Let now $s\in [3,p-4]$. Then from $x_1x\in D$ it follows that $\{x_{2},x_{3},\ldots , x_{s-1} \}\rightarrow x$. Together with statement \textbf{(i)} this implies that  
$$
\{x_{s-2},x_{s-1}\}\rightarrow x\rightarrow \{x_{s+1},x_{s+2}\}.
$$
 By  statements \textbf{(f)} and \textbf{(g)} we have $x_{s+1}x_s$,\, $x_sx_{s-1}\in D$,\, $x_sx_{s+2}\notin D$ and $x_{s-2}x_s\notin D$. 

It is not difficult to show that
$$
 A(x_s,\{x_{s-2},x_{s+2}\})=\emptyset. \eqno (29)
$$

Indeed, if $x_sx_{s-2}\in D$, then, since $\{x_{s-2},x_{s-1}\}\rightarrow x$, for the path $x_{s+1}x_{s+2}\ldots x_pyx_1x_2\ldots$ $ x_{s-1}$ and for the vertex $x_s$  we will have the considered Case 1, and if $x_{s+2}x_s\in D$, then in digraph $\overleftarrow D$ for the path $x_{s-1}x_{s-2}\ldots  x_1yx_px_{p-1}\ldots x_{s+2}x_{s+1}$ and for the vertex $x_s$ again we will have the considered Case 1 ($k\leq p-2$) and (29) holds.

 Recall that $A(x_s,x_p)=\emptyset$ by statement \textbf{(b)}. Together with (29) and $A(x_s,x)=\emptyset$ this implies that 
$$A(x_s,\{x,x_p,x_{s-2},x_{s+2}\}=\emptyset,$$ 
i.e., the vertex $x_s$ is not adjacent with at least 4 vertices of cycle $x_1x_2\ldots  x_{s-1}xx_{s+1}\ldots x_pyx_1$, this is contrary to statement \textbf{(i)}. The proof of Claim 10 is completed. \fbox \\\\

  Similarly to Claim 10 ($x_{p-1}y\notin D$), we can show the following claim:\\

\noindent\textbf{Claim 11}. $x_{p-1}x\notin (G)$, $xx_2\notin D$  and  $yx_2\notin D$. \fbox \\\\

Now let us complete the proof of Theorem 2. Without loss of generality, we may assume that $od(x)=n$. It follows that $x\rightarrow \{x_i,x_{i+1}\}$ for some $i\in [1,p-1]$. Using Claims 10, 11 and 1(ii) we see that $i\geq 3$ and  $x_1\rightarrow \{x,y\}$. Therefore  $A(x,x_l)=\emptyset$,\, $x\rightarrow \{x_{l+1},x_{l+2}\}$ and $x_{l-1}x\in D$ for some $l\in [2,i-1]$. So, for the path $x_{l+1}x_{l+2}\ldots x_pyx_1\ldots x_{l-1}$ we have $A(x,x_l)=\emptyset$, \,  $x\rightarrow \{x_{l+1},x_{l+2}\}$ and $x_{l-1}x_l$,  $x_lx_{l+1}$, $x_{l-1}x\in D$. This is a contradiction to  Claim 11 ($xx_2\notin D$) that $xx_{l+2}\notin D$. The proof of  Theorem 2 is completed. \fbox \\\\

\noindent\textbf {4.  Cycles of length  3 and 4 in digraph $D$.}\\

The next two results will be used in the proof of Theorem 3.\\

\noindent\textbf {Theorem A} (R. H\"{a}ggkvist, R. J. Faudree, R.H. Schelp [20]). Let $G$ be an undirected graph on $2n+1\geq 7$ vertices with minimum degree  at least $n$. Then precisely one of the following hold: (i) $G$ is pancyclic;   \, (ii) $ G\equiv (K_n\cup K_n)+K_1$; or (iii) $ K_{n,n+1}\subseteq G \subseteq K_n+\overline K_{n+1}$.\fbox \\\\

\noindent\textbf {Theorem B} (C. Tomassen [30]). Let $D$ be a strongly connected digraph on $p\geq 3$  vertices. If for each pair $x,y$ of nonadjacent  distinct vertices $d(x)+d(y)\geq 2p$, then $D$ is pancyclic or $p$ is even and $D\equiv K_{p/2,p/2}^*$.\fbox \\\\

Now we difine the digraphs $C^{*}_{6}(1)$, $H^{'}_{6}$ and  $H^{''}_{6}$ as folllows:

(i) $V(C^{*}_{6}(1))= \{x_1,x_2, \ldots , x_6 \} $  and  $A(C^{*}_{6}(1))= \{x_ix_{i+1}, x_{i+1}x_i / i\in [1,5] \}\cup \{ x_1x_6, x_6x_1, x_1x_3, x_1x_5,$ $ x_2x_4, x_6x_4 \}$;

(ii) $V(H^{'}_{6})=V(H^{''}_{6})=\{x,y,z, u,v,$ $w \}$,  $A(H^{'}_{6})$ $=\{ux,xu,xv, vx, yz, zy,zw,wz,xw,xy,uz,vz,wu,$ $wv,yu,yv \}$ and $A(H^{''}_{6})=\{ux,xu,xw,$ $xy, vx, vz,vw,wv,wu,zw,zy,yz,uz,yu,yv \}$.\\  
 
\noindent\textbf {Theorem 3}. Let $D$ be a digraph on $p\geq 5$  vertices with  minimum degree  at least $p-1$ and with  minimum semi-degree   at least $p/2-1$. Then the following hold: 

(i) $D$ contains a cycle of length 3 or  $p=2n$ and $D\subseteq K_{n,n}^*$  or else $D\in \{ C_5^*, K_{n,n+1}^* \}$;

(ii) $D$ contains a cycle of length 4  or $D\in \{ C_5^*, H^{'}_{6}, H^{''}_{6}, C^{*}_{6}(1), H(3,3), [(K_2 \cup K_2)+K_1]^*  \}$ .
 
\noindent\textbf {Proof}. Using Theorems A and B, we see that Theorem 3 is true if $D$ is a symmetric digraph. Suppose that $D$ is not symmetric digraph. If $D$ contains no cycle of length 3, then it is not difficult to show that $p=2n$ and $D\subseteq K_{n,n}^*$ (we leave the details to the reader).

 Assume that $D$ contains no cycle of length 4. For each arc $xy\in D$ put 
$$
S(x,y):=I(x)\cap O(y) \quad \hbox {and} \quad  E(x,y):=V(D)\setminus (O(y)\cup I(x)\cup \{x,y\}).
$$

Since $D$ has no  cycle of length 4, we see that
$$
A(O(y)\setminus \{x\} \rightarrow I(x)\setminus \{y\})=\emptyset. \eqno {(30)}
$$

 Let us consider the following cases.\\

\noindent\textbf {Case 1}. There is an arc $xy\in D$  such that $yx\notin D$ and $od(y)\geq n$ or $id(x)\geq n$, where $n:=\lfloor p/2 \rfloor$.

Without loss of generality, we can assume that $od(y)\geq n$  (if $id(x)\geq n$, then we will consider the digraph $\overleftarrow D$). Then from (30) and Lemma 4(ii) it follows that 
$$
I(x)\subseteq O(y), \quad I(x)=S(x,y) \quad \hbox {and} \quad A(\langle S(x,y)\rangle)=\emptyset. \eqno {(31)}
$$

 We now  shall prove that
$$
I(x)=O(y). \eqno {(32)}
$$
\noindent\textbf {Proof of (32)}. Assume that (32) is not true. Then  $ O(y)\setminus I(x)\not= \emptyset$ by (31), and let $z\in O(y)\setminus I(x)$. By (30),  $A(z\rightarrow \{x\}\cup I(x))=\emptyset$. From this  and Lemma 4(ii) it follows that  $|\{x\}\cup I(x)| =p/2$, $p=2n\geq 6$, $id(x)=n-1$ and  
$$
z\rightarrow V(D)\setminus (\{x,z\}\cup I(x)),   \eqno {(33)}
$$
in particular, $zy\in D$, $od(z)=n-1$ and $id(z)\geq n$. If $xz\in D$, then $C_4=xzyux$, where $u\in S(x,y)$, contradicting the our assumption. Therefore $xz\notin D$. From $id(x)=n-1$, $id(z)\geq n$ and Lemma 4(ii) it follows that $uz\in D$ for some vertex $u\in I(x)$. Now it is not difficult to see that $O(y)\setminus I(x)=\{z \}$, $E(x,y)\not=\emptyset$ and $A(E(x,y)\rightarrow y)=\emptyset$.

Suppose first that for each vertex $v\in I(x)$ there is a vertex $v_1\in E(x,y)$ such that $v_1v\in D$. Hence  $A(I(x)\rightarrow y)=\emptyset$ by (33) and $C_4\not\subset D$. Therefore  $A(I(x)\cup E(x,y)\rightarrow y)=\emptyset$, $|E(x,y)|=1$ and $n=3$. Let $E(x,y):=\{w\}$ and $I(x):=\{u,v\}$. Note that $w \rightarrow \{u,v\}$. Now it is not difficult to see  that if $wz\in D$, then $D \equiv H^{'}_{6}$ and if $wz\notin D$, then $D \equiv H^{''}_{6}$.  

Suppose next that $A(E(x,y)\rightarrow v)= \emptyset$ for some $v\in I(x)$. Then $|E(x,y)|=1$ by (31) and $n=3$, $xv\in D$, $vy\notin D$ and $v \rightarrow \{z, w\}$, where $w\in E(x,y)$. Now it is easy to see that $O(w)=\{z \}$. Therefore $od(w)\leq 1$, a contradiction.
 This proves (32), i.e., $I(x)=O(y)=S(x,y)$.\\

\noindent\textbf {Subcase 1.1}. $ A(x\rightarrow S(x,y))\not= \emptyset$.

Let $xu\in A(x\rightarrow S(x,y))$. If $uy\in D$, then  $C_4=xuyu_1x$, where $u_1\in S(x,y)\setminus \{u\}$,  a contradiction. So we may assume that $ uy \notin D$. From $od(y)\geq n$, (32) and (31) we get that $u\rightarrow E(x,y)$ and $od(y)=n$. It is not difficult to see  that $ A(E(x,y)\rightarrow (S(x,y)\setminus \{u\}))= \emptyset$ (otherwise $C_4\subset D$). Then  $E(x,y):=\{w\}$, $wu, wy\in D$ since $wx\notin D$, and $xv\in D$, where $v\in S(x,y)\setminus \{u\}$. Then $vy\in D$ or $vw\in D$. In both case we obtain a cycle of length 4, which is a contradiction. \\

\noindent\textbf {Subcase 1.2}. $ A(x\rightarrow S(x,y))= \emptyset$.

We can assume that $A(S(x,y)\rightarrow y)=\emptyset$ (otherwise in  $\overleftarrow D$ we will have Subcase 1.1). From $od(y)\geq n $,  by (32) and Lemma 4(ii) we have $E(x,y)\not= \emptyset$ and $x\rightarrow E(x,y) \rightarrow y$. Therefore $C_4=xzyux$, where $z\in E(x,y)$ and $u\in S(x,y)$, which is a contradiction and completes the discussion of Case 1.\\

\noindent\textbf {Case 2}. For each arc $xy\in D$ if $yx\notin D$,  then $od(y)<n$ and $id(x)<n$.

  From conditions of theorem it follows  easily that $od(y)=id(x)=n-1$ and  $p=2n\geq 6$. If $S(x,y)= \emptyset$, then using (30) and Lemma 4 (ii) it is easy to see that  $C_4\subset D$ or $D\in H(3,3)$. Assume that
 $S(x,y)\not= \emptyset$. Since $id(y)$ and $od(x)\geq n$, we can assume that $O(y)\rightarrow y$ and $x\rightarrow I(x)$ (otherwise for some arc $ux$ or $yv$ we have the considered Case 1). Hence it is easy to see that $|S(x,y)| =1$, $I(x)\not= O(y)$  and  $|E(x,y)| =1$. Let $E(x,y):=\{w\}$ and $S(x,y):=\{z\}$. From (30) it follows that $O(y)\setminus \{z \} \rightarrow w \rightarrow I(x)\setminus \{z\}$. From this, we obtain  $A(z,w)=\emptyset$ since $C_4\not\subset D$. Now it is not difficult to see that for some $v\in I(x)-\{z\}$ (or $u\in O(y)-\{z\}$)  $vz\in D$ (or $zu\in D$). Without loss of generality we may assume that $zu\in D$. From this we have $A(I(x)\setminus \{z\}\rightarrow  z)=\emptyset$, $wy\notin D$ and $n=3$. Hence $wu, vw, xw, vu \in D$ and $D\equiv  C^{*}_{6}(1)$. This completes the proof of  Theorem 3. \fbox \\\\

In [13], we proved the following:\\

\noindent\textbf{Theorem.} Let $D$ be a digraph on $ p\ge10$ vertices  with  minimum degree  at least $p-1$ and with minimum semi-degree at least $p/2-1$ ($n:=\lfloor p/2 \rfloor$). Then  $D$ is pancyclic unless 
 
$$ p=2n+1 \quad \hbox {and} \quad K_{n,n+1}^*\subseteq D \subseteq (K_n+\overline K_{n+1})^* \quad \hbox {or} \quad p=2n \quad \hbox {and} \quad G\subseteq K^*_{n,n}$$ 
 or else
$$
D\in H(n,n)\cup H(n,n-1,1)\cup \{ [(K_n \cup K_n)+K_1]^*, H(2n), H^ \prime (2n)\}.
$$
\\

\noindent\textbf {References}\\

[1] J. Bang-Jensen, G. Gutin, Digraphs: Theory, Algorithms and Applications, Springer, 2000.

[2] J. Bang-Jensen, G. Gutin, H. Li, Sufficient conditions for a digraph to be hamiltonian, J. Graph Theory 22 (2) (1996) 181-187.

[3] J. Bang-Jensen, Y. Guo, A. Yeo, A new sufficient condition for a digraph to be Hamiltonian, Discrete Applied Math. 95 (1999) 61-72.

[4] J. Bang-Jensen, Y. Guo, A note on vertex pancyclic oriented graphs, J. Graph Theory 31 (1999) 313-318.

[5] A. Benhocine, Pancyclism and Meyniel's conditions, Discrete Math. 58 (1986) 113-120.

[6] J. A. Bondy, C. Thomassen, A short proof of Meyniel's theorem, Discrete Math. 19 (1977) 195-197.

[7]  D. Christofides, P. Keevash, D. K\"{u}hn, D. Osthus, A semi-exact degree condition for Hamilton cycles in digraphs, submitted for publication.

[8] S. Kh. Darbinyan, Pancyclic and panconnected digraphs, Ph. D. Thesis, Institute  Mathematici Akad. Navuk BSSR, Minsk, 1981 (in Russian). 

[9] S. Kh. Darbinyan, Cycles of any length in digraphs with large semidegrees, Akad. Nauk Armyan. SSR Dokl. 75 (4) (1982) 147-152 (in Russian).

[10] S. Kh. Darbinyan, Pancyclicity of digraphs with the Meyniel condition, Studia Sci. Math. Hungar., 20 (1-4) (1985) 95-117 (in Russian).

[11] S. Kh. Darbinyan, Pancyclicity of digraphs with large semidegrees, Akad. Nauk Armyan. SSR Dokl. 80 (2) (1985) 51-54 (see also in  Math. Problems in Computer Science 14 (1985) 55-74) (in Russian).

[12] S. Kh. Darbinyan, A sufficient condition for the Hamiltonian property of digraphs with  large semidegrees, Akad. Nauk Armyan. SSR Dokl. 82 (1) (1986) 6-8 (in Russian).

[13] S. Kh. Darbinyan, On the  pancyclicity of digraphs with large semidegrees,  Akad. Nauk Armyan. SSR Dokl. 83 (3) (1986) 99-101 (in Russian).

[14] S. Kh. Darbinyan, A sufficient condition for  digraphs to be  Hamiltonian, Akad. Nauk Armyan. SSR Dokl. 91 (2) (1990) 57-59 (in Russian).

[15] S. Kh. Darbinyan, I. A. Karapetyan, On vertex pancyclicity oriented graphs, CSIT Conference, Yerevan, Armenia (2005) 154-155 (see also in  Math. Problems in Computer Science, 29  (2007) 66-84) (in Russian).

[16] S. Kh. Darbinyan, I. A. Karapetyan, On the large cycles through any given vertex in oriented graphs, CSIT Conference, Yerevan, Armenia (2007) 77-78 ( see also in  Math. Problems in Computer Science, 31 (2008) 90-107 (in Russian)).

[17] S. Kh. Darbinyan, K. M. Mosesyan, On  pancyclic regular oriented graphs, Akad. Nauk Armyan SSR Dokl. 67 (4) (1978) 208-211 (in Russian).

[18] A. Ghouila-Houri, Une condition suffisante d'existence d'un circuit hamiltonien, C. R. Acad. Sci. Paris Ser. A-B 251 (1960) 495-497.

[19] G. Gutin, Characterizations of vertex pancyclic and pancyclic ordinary complete multipartie digraphs,  Discrete Math., 141 (1-3) (1995) 153-162.

[20] R. H\"{a}ggkvist, R. J. Faudree, R. H. Schelp, Pancyclic graphs-connected Ramsey number, Ars Combinatoria 11 (1981) 37-49.

[21] R. H\"{a}ggkvist, C. Thomassen, On pancyclic digraphs, J. Combin. Theory Ser. B 20 (1976) 20-40.

[22] B. Jackson,  Long paths and cycles in oriented graphs, J. Graph Theory 5 (2) (1981) 145-157.

[23] P. Keevash, D. K\"{u}hn, D. Osthus, An exact minimum degree condition for Hamilton cycles in oriented graphs, J. London Math. Soc. 79 (2009) 144-166. 

[24] L. Kelly, D. K\"{u}hn, D. Osthus, Cycles of given length in oriented graphs, J. Combin. Theory Ser. B 100 (2010) 251-264.

[25]  D. K\"{u}hn, D. Osthus, A. Treglown, Hamiltonan degree sequences in digraphs, J. Combin. Theory Ser. B 100 (2010) 367-380.

[26] M. Meyniel, Une condition suffisante d'existence d'un circuit hamiltonien dans un graphe oriente, J. Combin. Theory Ser. B 14 (1973) 137-147.

[27] C. St. J. A. Nash-Williams, Hamilton circuits in graphs and digraphs, in: The Many Facets of Graph Theory, Springer- Verlag Lecture Notes, vol. 110, Springer Verlag (1969) 237-243.

[28] M. Overbeck-Larisch, A theorem on pancyclic-oriented graphs, J. Combin. Theory Ser. B 23 (2-3) (1977) 168-173.

[29] Z. M. Song, Pancyclic oriented graphs, J. Graph Theory 18 (5) (1994) 461-468.

[30] C. Thomassen, An Ore-type condition implying a digraph to be pancyclic, Discrete Math. 19 (1) (1977) 85-92.

[31] C. Thomassen, Long cycles in digraphs,  Proc. London Math. Soc. (3) 42 (1981) 231-251.

[32] D. R. Woodall, Sufficient conditions for circuits in graphs, Proc. London Math. Soc. 24 (1972) 739-755.

[33]   C. Q. Zhang, Arc-disjoint circuits in digraphs, Discrete Math. 41 (1982) 79-96.\\

\end{document}